\documentclass{article}
\usepackage{spconf,amsmath,graphicx}
\usepackage{graphicx}
\usepackage{float}
\usepackage{amsfonts}
\usepackage{amsmath}
\usepackage{amssymb}
\usepackage{cite}
\usepackage{bm}
\usepackage{mathrsfs}
\usepackage{epsfig}
\usepackage{algorithm}
\usepackage{algorithmic}
\usepackage{amsthm}
\usepackage{color}

\title{Recovery of  Sparse  Matrices via Matrix Sketching}
\name{Thakshila Wimalajeewa$^\dagger$, Yonina C. Eldar$^*$ and Pramod K. Varshney$^\dagger$ \thanks{The work of T. Wimalajeewa and P. K. Varshney  are supported by the National Science
Foundation (NSF) under Grant No. 1307775. The work of Y. C.  Eldar is supported in part by the Intel Collaborative Research Institute for Computational Intelligence (ICRI-CI), by the Ollendorf Foundation, by the Israel Science Foundation under grant 170/10, and by the Semiconductor Research Corporation (SRC) through the Texas Analog Center of Excellence (TxACE).}}
\address{$^\dagger$Dept. of Electrical Eng. and Comp.  Science,
   Syracuse University,
      Syracuse NY, USA \\
      $^*$Dept. of Electrical Engineering,
   Technion-Israel Institute of Technology,
      Technion City, Haifa,  Israel}
\begin{document}
\ninept
\maketitle
\begin{abstract}
In this paper, we consider the problem of recovering an unknown  sparse matrix $\mathbf X$ from the matrix sketch   $\mathbf Y = \mathbf A \mathbf X \mathbf B^T$. The dimension of $\mathbf Y$ is less than that of $\mathbf X$, and $\mathbf A$ and $\mathbf B$ are known matrices.  This problem can be solved using standard compressive sensing (CS) theory  after converting it to  vector form using the Kronecker operation. In this case,  the measurement matrix assumes a Kronecker product structure. However,  as the matrix dimension increases the associated computational complexity  makes its use prohibitive.  We extend two algorithms, fast iterative shrinkage threshold algorithm (FISTA)  and orthogonal matching pursuit (OMP) to solve this  problem in  matrix form  without employing  the Kronecker product. While both FISTA and OMP with matrix inputs are shown to be equivalent in performance to their vector counterparts with the Kronecker product, solving them in matrix form is shown to be computationally more efficient. We show that the computational gain achieved by FISTA with matrix inputs over  its vector form is more significant compared to that  achieved by OMP.  \end{abstract}
\begin{keywords}
Compressive sensing, Sparse matrix recovery, $l_l$ norm minimization, FISTA, OMP
\end{keywords}

\section{Introduction}
We consider the problem of recovering an unknown matrix $\mathbf X$ from the following observation model
\begin{eqnarray}
\mathbf Y = \mathbf A \mathbf X \mathbf B^T\label{obs_1}
\end{eqnarray}
where $\mathbf X\in \mathbb R^{N\times N}$, $\mathbf A \in \mathbb R^{M\times N}$, $\mathbf B \in \mathbb R^{L\times N}$ and $\mathbf A^T$ denotes the transpose of the matrix $\mathbf A$.
This problem  has been studied by many researchers in different contexts for arbitrary matrices  $\mathbf X$ \cite{Penrose1,Dai2,Peng1}.

In many applications dealing with high dimensional data, \emph{sparsity} is one of the low dimensional structures widely observed.
  Most popular transforms applied to 2-dimensional signals are in  the form of (\ref{obs_1}) where compression is obtained by a transformation of rows followed by a transformation of columns of the data matrix \cite{Caiafa1,Fang2,Rivenson1,Dasarathy1}.
With an arbitrarily  distributed sparse matrix $\mathbf X$ in which each column/row has only a few non zeros, a natural question to ask is whether it is possible to design  sensing matrices in the form of (\ref{obs_1}) so that $\mathbf X$ can  be uniquely  recovered from  $\mathbf Y$ when $M,L < N$.  Sparse signal recovery has  attracted much attention in the recent literature in the context of \emph{compressive sensing (CS)} \cite{candes1,donoho1,Eldar_B1}.   In the standard CS framework, a commonly used mechanism is to stack the high dimensional data into  vector form to recover the sparse vector uniquely from an  underdetermined linear system \cite{candes1,donoho1}.

The observation model (\ref{obs_1}) can be equivalently  written in  vector form using
 Kronecker  products as:
\begin{eqnarray}
\mathbf y = \mathbf C  \mathbf x\label{obs_2}
\end{eqnarray}
where $\mathbf y = \mathrm{vec}(\mathbf Y) \in \mathbb R^{ML}$, $\mathbf C = \mathbf B \otimes \mathbf A  \in \mathbb R^{ML\times N^2} $, $\mathbf x = \mathrm{vec}(\mathbf X) \in \mathbb R^{N^2}$, $\otimes$ denotes the Kronecker operator and $\mathrm{vec}(\mathbf X)$ is a column vector that vectorizes the matrix $\mathbf X$ (i.e. columns of $\mathbf X$ are stacked one after the other).
The sensing matrix in (\ref{obs_2}) has a special structure, i.e.,  it can be represented as a Kronecker product of two matrices $\mathbf A$ and $\mathbf B$.
It has been shown \cite{Duarte4,Jokar1,Duarte5,Jokar2} that the sparse signal $\mathbf x$ from (\ref{obs_2}) can be recovered  by solving the  following  $l_1$ norm minimization problem
\begin{eqnarray}
\min ||\mathbf x||_1~ s.t.~ \mathbf C \mathbf x = \mathbf y \label{l_1_norm_min}
\end{eqnarray}
 under certain conditions on the matrices  $\mathbf A$ and $\mathbf B$ where $||\mathbf x||_p$    denotes the  $l_p$ norm of $\mathbf x$. In particular, these results imply that the capability of recovering $\mathbf x$ based on (\ref{obs_2}) is ultimately determined by the worst behavior of $\mathbf A$ or $\mathbf B$. Also,  this approach is computationally complex  especially when the matrix dimension $N$  increases\cite{Rivenson1,Dasarathy1}.

  Several recent papers addressed the problem of recovering a sparse $\mathbf X$ from (\ref{obs_1}) without employing the Kronecker product. In \cite{Dasarathy1},  it was shown that a unique solution for $\mathbf X$ can be found when $\mathbf X$ is distributed  sparse under certain conditions on $\mathbf A$ and $\mathbf B$   by solving the following optimization problem:
\begin{eqnarray}
\min ||\mathbf X||_1 ~ \mathrm{s}.  ~\mathrm{t}. ~ \mathbf A \mathbf X \mathbf B^T = \mathbf Y\label{matrix_l1}
\end{eqnarray}
where $||\mathbf X||_1$ is  the $l_1$ norm of  $\mathrm{vec}(\mathbf X)$. The authors derive   recovery conditions when the matrices  $\mathbf A$ and $\mathbf B$ contain binary elements which are better than  that obtained via the Kronecker product approach. In \cite{Rivenson1}, the authors discuss  advantages in terms of
computational, storage, calibration and implementation while  solving (\ref{matrix_l1}) in matrix form compared to that with vector form. However, no specific  algorithm was developed to solve for $\mathbf X$.  In \cite{Fang2}, a version of orthogonal matching pursuit (OMP) (dubbed 2D OMP) is presented  to find a sparse $\mathbf X$ in the matrix form (\ref{obs_1}) when $\mathbf A= \mathbf B$.

%
%
%

Our goal in this paper is to develop algorithms to solve for sparse $\mathbf X$ from (\ref{obs_1}) without the employment of Kronecker products.  We extend fast iterative shrinkage threshold algorithm (FISTA) \cite{Beck1,Yang1}  developed for the vector case to sparse matrix recovery with  matrix inputs.  We further consider a greedy based approach via OMP to find the sparse solution. We show that  both algorithms with matrix inputs are equivalent to their vector counterparts obtained via Kronecker products in terms of performance. However,   the computational complexity of the  matrix approach is shown to be much less, especially with FISTA, compared to solving the problem in vector form.


\section{Sparse Matrix Recovery via $\l_1$ Norm Minimization}\label{matrix_l1norm}
\subsection{Vector formulation}
While numerous algorithms have been proposed in the literature to solve (\ref{l_1_norm_min}), in this paper we consider FISTA as discussed in \cite{Beck1,Yang1}. We consider the noisy observation model so that FISTA with vector inputs as given in Algorithm \ref{algo_FISTA_vec} \cite{Yang1}, is the solution of
\begin{eqnarray}
 \underset {\mathbf x}{\min} \left\{\frac{1}{2} ||\mathbf y - \mathbf C \mathbf x||_2^2 + \lambda ||\mathbf x||_1\right\}\label{l1_norm}
\end{eqnarray}
where $\lambda$ is a regularization parameter. In Algorithm \ref{algo_FISTA_vec}, $L_f=||\mathbf C||_2$ is the Lipschitz  constant of $\nabla f(\mathbf x)$  where  $||\mathbf C||_2$ denotes the spectral norm of $\mathbf C$, $\nabla$ denotes the gradient operator,  and  $f(\mathbf x) = \frac{1}{2}||\mathbf y - \mathbf C \mathbf x||_2^2$, and
\begin{eqnarray}
\mathrm{soft}(\mathbf u, a) =
\begin{array}{ccc}
\mathrm{sgn}(\mathbf u_i)(|\mathbf u_i| - a)_+
\end{array}
\end{eqnarray}
 for  $i=1,\cdots,N^2$
where $\mathbf u_i$ is the $i$-th element of $\mathbf u$, $x_+$ equals $x$ if $x>0$ and equals $0$ otherwise.
\begin{algorithm}
\textbf{Input:} observation vector $\mathbf y$, measurement matrix $\mathbf C$\\
\textbf{output:} estimate for signal, $\hat{\mathbf x}$
\begin{enumerate}

\item Initialization:  $\mathbf x^{0} =\mathbf 0$, $\mathbf x^{1}=\mathbf 0$, $t_0=1$, $t_1=1$, $k=1$\\
    Initialize: $\lambda_1$, $\beta\in (0,1)$, $\bar\lambda > 0$
\item \textbf{while} not converged \textbf{do}
\item $\mathbf z^{k} = \mathbf x^k + \frac{t_{k-1}-1}{t_k} (\mathbf x^k - \mathbf x^{k-1})$
 \item $ \mathbf u^k = \mathbf z^k - \frac{1}{L_f}\mathbf C^T (\mathbf C \mathbf z^k - \mathbf y) $
     \item $\mathbf x^{k+1} = \mathrm{soft} \left(\mathbf u^k, \frac{\lambda_k}{L_f}\right)$
     \item $t_{k+1} = \frac{1+\sqrt{4t_k^2 + 1}}{2}$
     \item $\lambda_{k+1} = \max(\beta \lambda_k, \bar\lambda)$
     \item $k=k+1$
     \item \textbf{end while}
 \end{enumerate}
 $\hat{\mathbf x} = \mathbf x^{k} $
 \caption{FISTA for sparse signal recovery with vector inputs}\label{algo_FISTA_vec}
 \end{algorithm}

 The computational complexity of FISTA is dominated by step 4 in Algorithm \ref{algo_FISTA_vec}.  The matrix-vector multiplications require $\mathcal O(N^4ML +N^4+N^2ML)$ computations. Since $M,L \leq N$, the  complexity is in  the order of $\mathcal O(N^4ML)$.  Thus, FISTA  is feasible only when $N, M, L$ are fairly  small.

 \subsection{Matrix formulation}
With the noisy version of (\ref{obs_1}), we aim to solve the following $l_1$ norm minimization problem:
\begin{eqnarray}
 \underset{\mathbf X}{\min} ~ \left\{\frac{1}{2} ||\mathbf Y - \mathbf A \mathbf X \mathbf B^T||_F^2 + \lambda ||\mathbf X||_1 \right\} \label{FISTA_matrix}
\end{eqnarray}
where $\lambda$ is a regularization parameter and $||\mathbf A||_F$ is the Frobenius norm of $\mathbf A$.
In generalizing FISTA to solve (\ref{FISTA_matrix}), we follow a similar approach as discussed in \cite{Toh1}.
Consider the more general unconstrained optimization problem:
\begin{eqnarray}
\underset{\mathbf X\in \mathbb R^{N\times N}}{\min} ~ F(\mathbf X) + \lambda G(\mathbf X)
\end{eqnarray}
where $G(\cdot)$ is a proper, convex, lower semicontinuous function, and $F(\cdot)$ is a convex smooth (continuously differentiable) function on an open subset of $\mathbb R^{N\times N}$ containing $\mathrm{dom} G = \{\mathbf X | G(\mathbf X) < \infty\}$. We assume that $\mathrm{dom} G$  is closed and $\nabla F$ is Lipschitz continuous on $\mathrm{dom} G$ with Lipschitz constant $L_f$; i.e.
\begin{eqnarray}
||\nabla F(\mathbf X) - \nabla F(\mathbf Z) ||_F \leq L_f ||\mathbf X - \mathbf Z||_F, ~X, Z\in \mathrm{dom}~G. \label{lips_mat}
\end{eqnarray}
When $F(\mathbf X) = \frac{1}{2} ||\mathbf Y - \mathbf A \mathbf X \mathbf B^T||_F^2$, it can be shown that $||\nabla F(\mathbf X) - \nabla F(\mathbf Z) ||_F = ||\nabla f(\mathbf x) - \nabla f(\mathbf z) ||_2$ and $||\mathbf X - \mathbf Z||_F = ||\mathbf x - \mathbf z||_2$ where $\mathbf z = \mathrm{vec}(\mathbf Z)$,   $f(\mathbf x) = \frac{1}{2}||\mathbf y - \mathbf C \mathbf x||_2^2$  and $\mathbf C=\mathbf B \otimes \mathbf A$  are as defined before.  Thus, the Lipschitz constant of  $\nabla F(\mathbf X)$ is the same as  $\nabla f(\mathbf x)$, and we use the same notation $L_f$ as used in Algorithm \ref{algo_FISTA_vec}.

Consider the following  quadratic approximation of $F(\cdot)$ at $\mathbf Z$ for any $\mathbf Z \in \mathrm{dom} G$:
\begin{eqnarray}
Q_L(\mathbf X, \mathbf Z) &:=& F(\mathbf Z) + \mathrm{tr}(\nabla F(\mathbf Z)^T (\mathbf X - \mathbf Z))\nonumber\\
& +& \frac{L_f}{2} ||\mathbf X - \mathbf Z||_F^2 + \lambda G (\mathbf X)\label{Q_L}
\end{eqnarray}
where $\mathrm{tr}(\cdot)$ denotes the trace of a matrix. We can rewrite  (\ref{Q_L})  as,
\begin{eqnarray}
Q_L(\mathbf X, \mathbf Z) =    \frac{L_f}{2} ||\mathbf X - \mathbf U(\mathbf Z)||_F^2 +  \lambda ||\mathbf X||_1 + \tilde F(\mathbf Z)\label{obj_matrix_Fista_mod}
\end{eqnarray}
where $\tilde F(\mathbf Z)$ is a function of only $\mathbf Z$ and
\begin{eqnarray}
\mathbf U(\mathbf Z) = \mathbf Z - \frac{1}{L_f}\nabla F(\mathbf Z)=\mathbf Z - \frac{1}{L_f} \mathbf A^T (\mathbf A \mathbf X \mathbf B^T - \mathbf Y)\mathbf B.
\end{eqnarray} Thus,
\begin{eqnarray*}
\underset{\mathbf X}{\arg\min}~ Q_L(\mathbf X, \mathbf Z) = \underset{\mathbf X}{\arg\min}  \left\{ \frac{L_f}{2} ||\mathbf X - \mathbf U(\mathbf Z)||_F^2 +  \lambda ||\mathbf X||_1 \right\}.
\end{eqnarray*}
Since both terms are element wise separable, we have
\begin{eqnarray}
\underset{\mathbf X\in \mathrm{dom} G}{\arg\min} ~ Q(\mathbf X, \mathbf Z)= \mathrm{soft}\left(\mathbf U(\mathbf Z), \frac{\lambda}{L_f}\right)\label{P_l_z}
\end{eqnarray}
where $\mathrm{soft}\left(\mathbf U(\mathbf Z), \frac{\lambda}{L_f}\right)$ denotes an element wise operation with  \begin{eqnarray*}
\mathrm{soft}\left(\mathbf W, L_0\right)=
\begin{array}{ccc}
\mathrm{sgn}(\mathbf W_{ij})(|\mathbf W_{ij}| - L_0)_+
\end{array}
\end{eqnarray*}
for all indices $i,j$ of the $N\times N$ matrix $\mathbf W$. These steps lead to a generalization of FISTA with matrix inputs, as given in Algorithm \ref{algo_FISTA_matrix}.

\begin{algorithm}
\textbf{Input:} observation matrix $\mathbf Y$, measurement matrices $\mathbf A$ and $\mathbf B$ \\
\textbf{Output:} estimate for sparse signal matrix, $\hat{\mathbf X}$
\begin{enumerate}
\item Initialization:  $\mathbf X^{0} =\mathbf 0$, $\mathbf X^{1}=\mathbf 0$, $t_0=1$, $t_1=1$, $k=1$\\
     Initialize: $\lambda_1$, $\beta\in (0,1)$, $\bar\lambda > 0$
\item \textbf{while} not converged \textbf{do}
\item $\mathbf Z^{k} = \mathbf X^k + \frac{t_{k-1}-1}{t_k} (\mathbf X^k - \mathbf X^{k-1})$
 \item $ \mathbf U^k = \mathbf Z^k - \frac{1}{L}\mathbf A^T (\mathbf A \mathbf Z^k\mathbf B^T - \mathbf Y) \mathbf B $
     \item $\mathbf X^{k+1} = \mathrm{soft} \left(\mathbf U^k, \frac{\lambda_k}{L_f}\right)$
     \item $t_{k+1} = \frac{1+\sqrt{4t_k^2 + 1}}{2}$
     \item $\lambda_{k+1} = \max(\beta \lambda_k, \bar\lambda)$
     \item $k=k+1$
     \item \textbf{end while}
 \end{enumerate}
 $\hat{\mathbf X} = \mathbf X^{k} $
 \caption{FISTA for sparse matrix  recovery with matrix inputs}\label{algo_FISTA_matrix}
 \end{algorithm}

 As in Algorithm \ref{algo_FISTA_vec}, the computational complexity is dominated by  step 4. The matrix-matrix multiplication at step 4 in Algorithm \ref{algo_FISTA_matrix} is performed with $\mathcal O(N^2(N+M+3L) + NML)$ computations. Since $M,L \leq N$, the worst case complexity  is $\mathcal O(N^3)$. Recall, that FISTA in vector form has worst case complexity of $\mathcal O(N^4ML)$. Thus, there is a $\mathcal O(NML)$ gain in the matrix version compared to the vector approach.

\subsection{Equivalence of Algorithms \ref{algo_FISTA_vec} and \ref{algo_FISTA_matrix}}
It is easy to see that $\mathbf z^k$ and $\mathbf x^{k+1}$ computed in steps 3 and 5 in Algorithm \ref{algo_FISTA_vec} are  the same as $\mathrm{vec}(\mathbf Z^k)$ and $\mathrm{vec } (\mathbf X^{k+1})$, respectively, if $\mathbf u^k = \mathrm{vec}(\mathbf U^k)$ where $\mathbf Z_k$, $\mathbf U_k$ and $\mathbf X^{k+1}$ are as computed at steps 3, 4 and 5 in Algorithm \ref{algo_FISTA_matrix} . Now,
\begin{eqnarray*}
&~&\mathrm{vec}(\mathbf A^T \mathbf A \mathbf Z^k \mathbf B^T \mathbf B - \mathbf A^T  \mathbf Y \mathbf B)\nonumber\\
&=&((\mathbf B^T \mathbf B)\otimes (\mathbf A^T \mathbf A)) \mathrm{vec}(\mathbf Z^k) - (\mathbf B^T\otimes \mathbf A^T )\mathrm{vec}(\mathbf Y)\nonumber\\
&=&(\mathbf B^T \otimes \mathbf A^T)(\mathbf B\otimes \mathbf A) \mathrm{vec}(\mathbf Z^k) - (\mathbf B^T\otimes \mathbf A^T )\mathrm{vec}(\mathbf Y)\nonumber\\
&=&(\mathbf C^T \mathbf C)\mathbf z^k -\mathbf C^T \mathbf y,
\end{eqnarray*}
where $\mathbf C= \mathbf B \otimes \mathbf A$.
Thus, $\mathbf u_k$ computed at step 4 in Algorithm \ref{algo_FISTA_vec} is the vectorized version of  $\mathbf U^k$ of  step 4 in Algorithm \ref{algo_FISTA_matrix}.
We conclude that  Algorithms  \ref{algo_FISTA_vec} and \ref{algo_FISTA_matrix} provide the same output, however, Algorithm \ref{algo_FISTA_matrix} is more efficient.

\begin{algorithm}
\textbf{Input:} observation matrix $\mathbf Y$, measurement matrices $\mathbf A$ and $\mathbf B$ \\
\textbf{Output:} index set $\Lambda$ containing locations of the non zero indices of the matrix $\mathbf X$, estimate for signal matrix $\hat{\mathbf X}$
\begin{enumerate}
\item Initialization:  residual $\mathbf R_0 =\mathbf Y$, index set   $\Lambda_0=\emptyset$, $t=1$
\item Find the two indices $\lambda_t=[\lambda_t(1)~ \lambda_t(2)]$ such that
    \begin{eqnarray}
    [\lambda_t(1)~\lambda_t(2)]& = &\underset{i,j}{\arg\max}~ |\mathbf b_j^T\mathbf R_{t-1}^T \mathbf a_i|\label{eq_OMP_matx}
    \end{eqnarray}
 \item Augment index set $\Lambda_t = \Lambda_t \cup \{\lambda_t\}$
     \item Find the new signal estimate
     \begin{eqnarray}
\mathbf x_t = \mathbf D_t^{-1} \mathbf d_t\label{x_t}
     \end{eqnarray}
     where $\mathbf D_t$ and $\mathbf d_t$ are as in (\ref{hat_t})
     \item Compute new  residual
     \begin{eqnarray}
     \mathbf R_t &=& \mathbf Y - \sum_{m=1}^{t} \mathbf x_t(m) \mathbf a_{\Lambda_t(m,1)} \mathbf b_{\Lambda_t(m,2)}^T\label{residual_M}
     \end{eqnarray}
     \item Increment $t$ and return to step $2$ if $t \leq d $, otherwise stop
   \item Estimated support set $\hat\Lambda=\Lambda_d$ \\
    Estimated signal matrix  $\hat{\mathbf X}$:   $(\Lambda_{d}(m,1), \Lambda_{d}(m,2))$-th component of $\hat{\mathbf X}$ is given by $\mathbf x_{d}(m)$ for $m=1=,\cdots, d$ while  rest of the elements are zeros.
 \end{enumerate}
 \caption{OMP with matrix inputs}\label{algo_OMP_matrix}
 \end{algorithm}

\section{Sparse Matrix Recovery via OMP}
Next, we consider the extension of standard OMP to  the matrix form    (\ref{obs_1}). We can write the observation $\mathbf Y$ in (\ref{obs_1}) as a summation of $N^2$ matrices as given below:
\begin{eqnarray}
\mathbf Y =\underset{i,j} {\sum} X_{ij} \mathbf a_i \mathbf b_j^T.\label{obs_sum}
\end{eqnarray}
When $\mathbf X$ is sparse with $d$ nonzeros, the summation in (\ref{obs_sum}) has only $d$ terms. Let $\Sigma_d$ denote the support of $\mathbf X$  so that $\mathbf X_{ij}$ is non zero for $i,j=1,\cdots, N$, and let $\bar\Sigma_d$ be its  complement.  We can write (\ref{obs_sum}) as $
\mathbf Y =\underset{(i,j) \in \Sigma_d} {\sum} X_{ij} \mathbf a_i \mathbf b_j^T.
$ Our goal is to recover $\mathbf X_{ij}$ for $(i,j)\in \Sigma_d$ in a greedy manner.
The proposed OMP version with matrix inputs is given in Algorithm \ref{algo_OMP_matrix}. In Algorithm \ref{algo_OMP_matrix}, $\Lambda_t$ contains estimated $(i,j)$ pairs up to $t$-th iteration in which the $m$-th pair is denoted by $(\Lambda_t(m,1), \Lambda_t(m,2))$  for $m=1,\cdots, t$. Once $\Lambda_t$ is updated as in step 3, the  signal is estimated solving the following optimization problem:
\begin{eqnarray}
\mathbf x_t = \underset{\mathbf x} {\arg\min}\parallel \mathbf Y - \sum_{m=1}^{t} x_m \mathbf a_{\Lambda_t(m,1)} \mathbf b_{\Lambda_t(m,2)}^T ) \parallel_F. \label{eq_frob_norm}
\end{eqnarray}
The solution of (\ref{eq_frob_norm}) is given by
\begin{eqnarray}
\mathbf x_t = \mathbf D_t^{-1} \mathbf d_t\label{hat_t}
\end{eqnarray}
where $\mathbf D_t$ is a $t\times  t$ matrix in which the $(m,r)$-th element is given by
 \begin{eqnarray}
 (\mathbf D_t)_{m,r} = \mathbf b_{\Lambda_t(r,2)}^T \mathbf b_{\Lambda_t(m,2)} \mathbf a_{\Lambda_t(m,1)}^T \mathbf a_{\Lambda_t(r,1)}
 \end{eqnarray}
 for $m,r=1,\cdots, t$ and
\begin{eqnarray}\mathbf d_t = [\mathbf b_{\Lambda_t(1,2)}^T \mathbf Y^T \mathbf a_{\Lambda_t(1,1)}~\cdots ~\mathbf b_{\Lambda_t(t,2)}^T \mathbf Y^T \mathbf a_{\Lambda_t(t,1)}]^T
\end{eqnarray} is a $t\times 1$ vector.
Then the new approximation at the $t$-th iteration  is given by
\begin{eqnarray}
  \mathbf Q_t &=& \sum_{m=1}^{t} \mathbf x_t(m) \mathbf a_{\Lambda_t(m,1)} \mathbf b_{\Lambda_t(m,2)}^T\label{Q_t}
  \end{eqnarray}
where $\mathbf x_t(m)$ denotes the $m$-th element of $\mathbf x_t$.

  Algorithm \ref{algo_OMP_matrix} is a trivial extension of the standard OMP (and was also considered in \cite{Fang2} for  $\mathbf A=\mathbf B$).
 \subsection{Computational complexity }
 As shown in \cite{Fang2} for $\mathbf A= \mathbf B$,  it can be easily verified that Algorithm \ref{algo_OMP_matrix} and the standard OMP \cite{tropp1}  with vector inputs (\ref{obs_2}) provide the same  performance at each iteration. However, the computational complexity of  Algorithm \ref{algo_OMP_matrix} is less than that of its vector counterpart. 
  Step 2 in Algorithm  \ref{algo_OMP_matrix} can be implemented as a matrix multiplication  $\mathbf A^T \mathbf R_{t-1} \mathbf B$. Thus, the computational complexity of this step is in the order of $\mathcal O(NML+N^2L)$. It is noted that, when implementing the standard OMP  as in \cite{tropp1} with vector form (\ref{obs_2}), the equivalent step is computed with  complexity of $\mathcal O(N^2ML)$. With respect to step 4 in Algorithm \ref{algo_OMP_matrix},  the matrix $\mathbf D_t$ requires $\mathcal O(t^2(M+L) )$ computations at the $t$-the iteration. The vector $\mathbf d_t$ requires  $\mathcal O(t(ML+M))$ computations.  Worst case complexity of the inverse operation is $\mathcal O(t^3)$. Matrix-vector multiplication in (\ref{x_t}) requires $\mathcal O(t^2)$ computations. Thus, at a given iteration,  worst case  complexity of step 4 in Algorithm \ref{algo_OMP_matrix} is in the order of $\mathcal O(tML)$. It can be shown that the worst case  computational complexity of the equivalent step of standard OMP with Kronecker products to estimate the signal at $t$-th iteration,  is in the order of $\mathcal O(t^2 ML)$.  Thus, steps 2 and 4 in Algorithm \ref{algo_OMP_matrix} provide  us with a computational gain over the equivalent steps of the standard OMP with Kronecker products.
Therefore, we conclude that Algorithm \ref{algo_OMP_matrix} is an efficient way to find sparse $\mathbf X$ from (\ref{obs_1})  compared to its vector counterpart (\ref{obs_2}) although both provide the same performance. It is further observed  that this  computational gain is not as  significant as with FISTA.

\section{Numerical Results}
In this section, we demonstrate the capability of recovering sparse $\mathbf X$ from observation model (\ref{obs_1}) via different algorithms and provide  insights into the computational gains achievable with  the matrix version.
First, we illustrate the performance of FISTA  with different choices for $\mathbf A$ and $\mathbf B$.
For numerical results, we assume that $\mathbf X$ is a distributed sparse matrix in which each column has a maximum of $K$ nonzeros and the locations are generated uniformly.  The values of nonzero entries are drawn from a uniform distribution in the range $ [-250, ~-200] \cup[200, ~250]$. We consider that  the  observation matrix $\mathbf Y$ in  (\ref{obs_1}) is corrupted by additive noise and the elements of the noise matrix  are assumed to be independently and identically distributed Gaussian random variables with mean zero and  variance  $\sigma_v^2$

\begin{figure}[htb]
 \centering
  \centerline{\includegraphics[width=9.0cm]{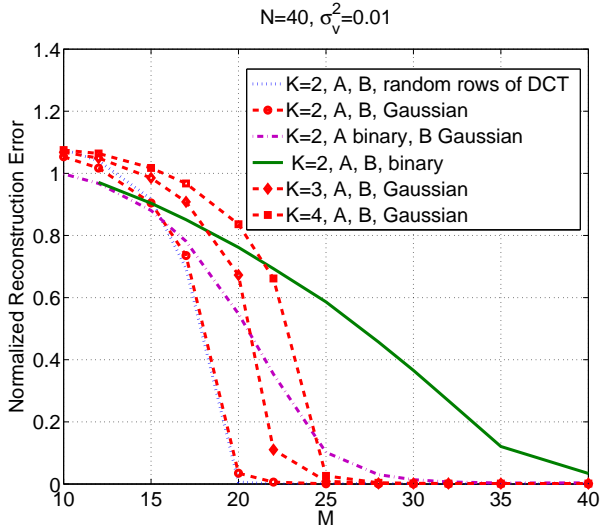}}
\caption{Normalized reconstruction error vs $M=L$ with FISTA and  different projection matrices, $N=40$, $\sigma_v^2=0.01$}
\label{fig_1}
\end{figure}
 In Fig.  \ref{fig_1},  we plot the normalized  reconstruction error $\frac{||\mathbf X - \hat{\mathbf X}||_F}{||\mathbf X||_F}$ vs $M$ when $M=L$ averaging over $500$ trials. We  let $N=40$, and $\sigma_v^2=0.01$.  In Fig. \ref{fig_1},  we illustrate two aspects. First,  for given $K$, different types of matrices $\mathbf A$ and $\mathbf B$ are examined.  We consider  independent random rows of the $N\times N$ DCT matrix, Gaussian, and binary matrices. In the case of a  Gaussian, elements are drawn from a normal ensemble and then orthogonalized. By binary matrix, we mean that the elements of the matrix can take values $1/N$ or $0$ with equal probability. Note that, random rows of DCT matrix and Gaussian matrix obey uniform uncertainty principle with good isometry constants in contrast to a binary matrix. When both matrices $\mathbf A$ and $\mathbf B$ are either random rows of DCT matrix or  zero mean Gaussian,  the recovery of the sparse matrix is guaranteed with less measurements compared to $N$.  When $\mathbf A$ and $\mathbf B$ are binary, the recovery is not so good, which is intuitive since binary matrices are not "good" compressive sensing matrices.  However, when $\mathbf A$ is binary and $\mathbf B$ is Gaussian, we see an improved performance compared to the case where both are binary.
 This  provides an insight that even when one matrix does not obey uniform uncertainty principle with good isometry constant, still the sparse matrix can be recovered reliably when the other matrix is a "good compressive sensing" matrix. We will further investigate this observation in our future work.
 Second, with a given type of matrices (in the case of $\mathbf A$, $\mathbf B$ are Gaussian)  we vary $K$. It is seen that, recovery capability of FISTA  does not degrade significantly as $K$ increases.

\begin{figure}[htb]
 \centering
  \centerline{\includegraphics[width=9.0cm]{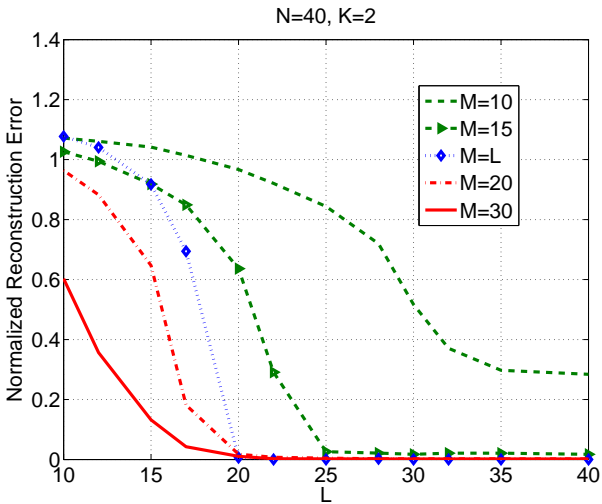}}
\caption{Normalized reconstruction error vs $L$ for given $M$ with FISTA, $\mathbf A$ and $\mathbf B$ contain independent random rows of the $N\times N$ DCT matrix,  $N=40$, $\sigma_v^2=0.01$}
\label{fig_01}
\end{figure}

In Fig. \ref{fig_01}, we plot the reconstruction error vs $L$ with FISTA keeping $M$ fixed. As a benchmark, the curve corresponding to $M=L$ is also plotted.
It can be seen that, when one dimension of the observation matrix  is fixed, an improved performance in terms of signal reconstruction error  is observed as the other dimension increases. However, when  $M$ is  very small (or below a certain value), complete recovery is not guaranteed even if $L=N$. This implies that when the dimension of the matrix $\mathbf A$ is fixed,  increasing the number of columns of $\mathbf B^T$ does not necessarily guarantee  complete recovery when $M$ is very small.

\begin{table}
\caption{Runtime (in $s$)  of FISTA with vector and matrix inputs}
\centering
\begin{small}
\begin{tabular}{|l|l|l|l|}
  \hline
 $~$   & N=20 & N=40 & N=60 \\
  \hline
  Matrix & $0.3863$ &  $1.3064$ & $6.5595$  \\
 Vector & $1.0987$ &  $26.4162$ & $142.7584$  \\
   \hline
 \end{tabular}
\label{table_comparison}
\end{small}
\end{table}
In Table \ref{table_comparison}, we compare the average runtime  with MATLAB (in Intel(R) Core(TM) i7-3770 CPU$@$ 3.40GHzz processor with 12 GB RAM) for FISTA for matrix and vector versions as the sparse matrix dimension $N$ varies given that the number of iterations in both Algorithms  \ref{algo_FISTA_vec}  and \ref{algo_FISTA_matrix} is fixed at the same value ($=10000$). We let $K=N/20$ and $M=L=N/2$. Matrices $\mathbf A$ and $\mathbf B$ are assumed to be Gaussian. It reflects the computational efficiency of the matrix approach compared to the vector approach   especially  as $N$ increases although both algorithms provide the same performance.

\begin{figure}[htb]
 \centering
  \centerline{\includegraphics[width=9.0cm]{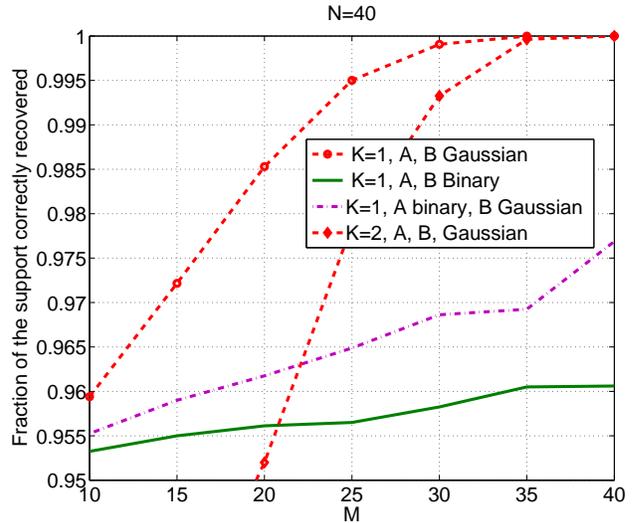}}
\caption{Fraction  of the support correctly recovered  vs $M=L$ with OMP  with different projection matrices with no noise}
\label{fig_3}
\end{figure}

To illustrate  the performance of OMP with matrix inputs, we plot the fraction  of the support correctly recovered with  Algorithm \ref{algo_OMP_matrix} for different choices for $\mathbf A$ and $\mathbf B$ with $K=1$ in Fig. \ref{fig_3}.  From Fig. \ref{fig_3}, it is again observed that, although not with the same scale as with FISTA,  the recovery capability can be improved when one projection matrix  is binary and the  other is  Gaussian compared to the case where both $\mathbf A$ and $\mathbf B$ are binary. Another observation is that even with a Gaussian matrix, as $K$ increases the performance of OMP degrades significantly leaving OMP not a better  choice when the sparsity level increases.

\section{Discussion}
In this paper, we showed numerically that recovering  $\mathbf X$ based on (\ref{obs_1}) in its matrix form is more computationally efficient than solving it after converting to vector form via Kronecker products  when $\mathbf X$ is arbitrarily distributed  sparse. We developed  matrix versions  of FISTA to solve $l_1$ norm minimization  in (\ref{matrix_l1}) efficiently and OMP to solve for $\mathbf X$ in a greedy manner. It has been shown that a  significant computational gain is achieved by FISTA with matrix form compared to its vector counterpart.  We further illustrated the recovery capability with different choices for projection operators.  The results provide insight into  the following. If a linear system of the form (\ref{obs_2}) can be  converted into a matrix form as in (\ref{obs_1}), the problem can be solved more efficiently without losing performance with respect to the original vector  form. Thus, it is worth investigating such scenarios where the matrix approach can be efficiently used to solve linear systems  which are computationally demanding otherwise.

%



\newpage
\bibliographystyle{IEEEtran}
\bibliography{IEEEabrv,bib1}

\end{document}